\documentclass[onefignum,final]{siamart171218}

\title{Least Angle Regression Coarsening in Bootstrap Algebraic Multigrid\thanks{This work was partially funded by Deutsche Forschungsgemeinschaft (DFG) Transregional Collaborative Research Centre 55 (SFB/TRR55)}}

\author{Karsten Kahl\thanks{Bergische Universiat\"at Wuppertal, Faculty of Mathematics and Natural Sciences, \texttt{\{kkahl,rottmann\}@math.uni-wuppertal.de}}\ \and Matthias Rottmann\footnotemark[2]}

\usepackage[utf8]{inputenc}
\usepackage[english]{babel}

\usepackage{amsmath,amsfonts}
\usepackage{mathabx}
\usepackage{textgreek}

\usepackage[vlined,boxruled,linesnumbered,algo2e,algosection]{algorithm2e}
\usepackage{xcolor}
\usepackage{tikz}
\usetikzlibrary{intersections,plotmarks}

\usepackage{csquotes}

\usepackage[
backend=biber, 
style=numeric,
defernumbers=true,
firstinits=true,
url=false,
doi=false,
isbn=false,
eprint=false,
maxbibnames=99,
]{biblatex}

\renewbibmacro{in:}{%
    \ifboolexpr{%
        test {\ifentrytype{article}}%
        or
        test {\ifentrytype{inproceedings}}%
    }{}{\printtext{\bibstring{in}\intitlepunct}}%
}

\newcommand{\SetAlgorithmStyle}{
  \setcounter{AlgoLine}{0}
  \SetKwData{Left}{left}\SetKwData{This}{this}\SetKwData{Up}{up}
  \SetKwInOut{Input}{Input}
  \SetKwInOut{Output}{Output}
  \ResetInOut{input}
  \SetKwComment{tcp}{//}{}
  \SetKwFor{For}{for}{}{end}
  \SetArgSty{}
  \DontPrintSemicolon
}

\usepackage{enumerate}

\newcommand{\linspace}[2][ ]{\mathbb{#2}^{#1}}
\newcommand{\innerprod}[3][2]{\langle #2,#3 \rangle_{#1}}
\newcommand{\norm}[2][2]{\| #2 \|_{#1}}

\newcommand{\ihat}{\mathit{\hat{\textit{\i}}}}
\newcommand{\icheck}{\mathit{\check{\textit{\i}}}}

\newcommand{\LARS}{least angle regression}
\newcommand{\CLARS}{Least angle regression}
\newcommand{\leastangle}{least angle}
\newcommand{\Leastangle}{Least angle}

\newcommand{\leastsquares}{least squares}

\newcommand{\LeastSquares}{Least Squares}

\crefalias{AlgoLine}{line}
\crefalias{algocf}{algorithm}

\makeatletter
\let\cref@old@stepcounter\stepcounter
\def\stepcounter#1{%
  \cref@old@stepcounter{#1}%
  \cref@constructprefix{#1}{\cref@result}%
  \@ifundefined{cref@#1@alias}%
    {\def\@tempa{#1}}%
    {\def\@tempa{\csname cref@#1@alias\endcsname}}%
  \protected@edef\cref@currentlabel{%
    [\@tempa][\arabic{#1}][\cref@result]%
    \csname p@#1\endcsname\csname the#1\endcsname}}
\makeatother

\makeatletter
\renewcommand{\algocf@caption@boxruled}{%
  \hrule
  \hbox to \hsize{%
    \vrule\hskip-0.4pt
    \vbox{   
       \vskip\interspacetitleboxruled%
       \unhbox\algocf@capbox\hfill
       \vskip\interspacetitleboxruled
       }%
     \hskip-0.4pt\vrule%
   }\nointerlineskip%
}%
\makeatother

\begin{document}

\maketitle

\begin{abstract}
The bootstrap algebraic multigrid framework allows for the adaptive construction of algebraic multigrid methods in situations where geometric multigrid methods are not known or not available at all. While there has been some work on adaptive coarsening in this framework in terms of algebraic distances, coarsening is the part of the adaptive bootstrap setup that is least developed. In this paper we try to close this gap by introducing an adaptive coarsening scheme that views interpolation as a local regression problem. In fact the bootstrap algebraic multigrid setup can be understood as a machine learning ansatz that learns the nature of smooth error by local regression. In order to turn this idea into a practical method we modify \leastsquares{} interpolation to both avoid overfitting of the data and to recover a sparse response that can be used to extract information about the coupling strength amongst variables like in classical algebraic multigrid. In order to improve the so-found coarse grid we propose a post-processing to ensure stability of the resulting \leastsquares{} interpolation operator. We conclude with numerical experiments that show the viability of the chosen approach.
\end{abstract}



\section{Introduction}
In this paper we establish a connection between the adaptive bootstrap algebraic multigrid setup, in particular its coarsening part, and data driven machine learning approaches. 

Algebraic multigrid methods have been introduced in~\cite{BranMcCoRuge1985,Bran1986,RugeStue1987,Stue1983} as a method to efficiently solve sparse linear systems of equations
\[
Ax = b,\ b\in\mathbb{C}^{n},\ A \in \mathbb{C}^{n \times n}
\] without the requirement of expert knowledge, e.g., the underlying physical model, the employed discretization scheme and/or geometry. Efficiency in algebraic multigrid methods is achieved by pairing a simple iterative scheme, the \textit{smoother}, with a coarse grid correction. Generically the error propagator of a two-grid algebraic multigrid method with Galerkin coarse grid construction can be written as
\begin{equation} \label{eq:eprop}
    E_{2g} = (I-MA) (I-P (P^H A P)^{-1} P^H A) (I-MA) \, ,
\end{equation} where it is the task of the algebraic multigrid setup to determine a suitable interpolation operator $P$. That is, one has to find suitable choices for the dimension of the coarse space, $n_{c}$, the sparsity pattern of $P\in\mathbb{C}^{n\times n_{c}}$ and its entries. Typically these tasks are split into two parts. Finding $n_{c}$ and the sparsity pattern of $P$ is often referred to as the \textit{coarsening} problem, while determining the entries of $P$ is known as the \textit{interpolation} problem.

In the classical algebraic multigrid approach~\cite{RugeStue1987} both problems are solved using the entries of $A$ and it has been shown that this is appropriate as long as $A$ has $M$-matrix structure, e.g., as a suitable discretization of an elliptic partial differential equation. The classical approach relies on quite restrictive assumptions on the underlying problem and therefore cannot be extended significantly beyond the M-matrix case. In recent years the scope of algebraic multigrid methods has been enlarged by the introduction of adaptivity~\cite{BranBranKahlLivs2011,BrezFalgMacLMantMcCoRuge2006,BrezFalgMacLMantMcCoRuge2004}. The fundamental idea of adaptive approaches is to guide the construction of the coarse space by either using spectral information on $A$ and/or the smoother, or simply using  the action of the smoother itself. While many of these approaches succeeded in addressing the interpolation problem, advances for the coarsening problem are scarce. 
Some approaches try to generalize the definition of strength of connection of the classical method~\cite{OlsoSchrTumi2010}. Others consider only binary relations of variables~\cite{BranBranKahlLivs2015b} and compute a strength of connection method in an adaptive fashion. This can also be said for recent algebraic aggregation approaches~\cite{BranChenKrauZika2013,LivnBran2012,NapoNota2016,Nota2010}, which mainly use binary relations as well.
Last, there are approaches based on compatible relaxation~\cite{BranFalg2010} that come closest to general applicability, but oftentimes do not mesh efficiently with the chosen approach for the definition of interpolation weights, i.e., the entries of $P$.

In this paper we propose a new way of solving the coarsening problem in the bootstrap algebraic multigrid framework. Based on the concept of \leastsquares{} interpolation we develop a \emph{\leastangle{} regression} coarsening scheme that can be fully integrated into the bootstrap framework and which utilizes only the information present in the small number of test vectors of the bootstrap process. To do so, we review the concept of \leastsquares{} interpolation in~\cref{sec:lsinterpolation} and show how it can be interpreted as a machine learning, i.e., regression, problem. We introduce an $\ell_{1}$ penalty term, also known as a \emph{lasso} term~\cite{Tibs1996}, into the \leastsquares{} interpolation and show how the modified problem can be solved efficiently by \emph{\leastangle{} regression}. 
We continue in~\cref{sec:larcoarsening} with the description of the overall coarsening strategy and show numerical results in~\cref{sec:numericalresults}.

\section{\LeastSquares{} Interpolation}\label{sec:lsinterpolation}
The bootstrap algebraic multigrid framework constructs a multigrid hierarchy by leveraging the information contained in a set of test vectors 
\[
\mathcal{V} = \{v^{(1)},\ldots,v^{(K)}\} \subset \mathbb{C}^n
\]
(cf.~\cite{BranBranKahlLivs2011}). The central part of the setup process is the calculation of interpolation weights by \leastsquares{} interpolation. Assuming that the set of variables $\Omega$ is split into a set of coarse variables $\mathcal{C}$ and fine variables $\mathcal{F} = \Omega \setminus \mathcal{C}$ and that the sparsity pattern of interpolation is known, i.e., each variable $i$ is equipped with a set of variables $\mathcal{C}_{i} \subset \mathcal{C}$ it interpolates from, the interpolation weights in \leastsquares{} interpolation are simply given by the weighted \leastsquares{} fit
\begin{equation}\label{eq:lsi}
    \sum_{k = 1}^{K} \omega_{k} \left(v_{i}^{(k)} - \sum_{j\in \mathcal{C}_{i}} p_{ij}v_{j}^{(k)}\right)^{2} \rightarrow \operatorname{min}.
\end{equation} Herein, $\omega_{k}$ is chosen to reflect the importance of test vector $v^{(k)}$. Interpolation for variables in $\mathcal{C}$ is defined by the identity. In accordance with~\cite{BranBranKahlLivs2011} we call the cardinality $|\mathcal{C}_{i}|$ of the set of interpolation variables the \emph{caliber} of interpolation.

\begin{figure}[htb]
\centerline{\includegraphics[width=1\textwidth]{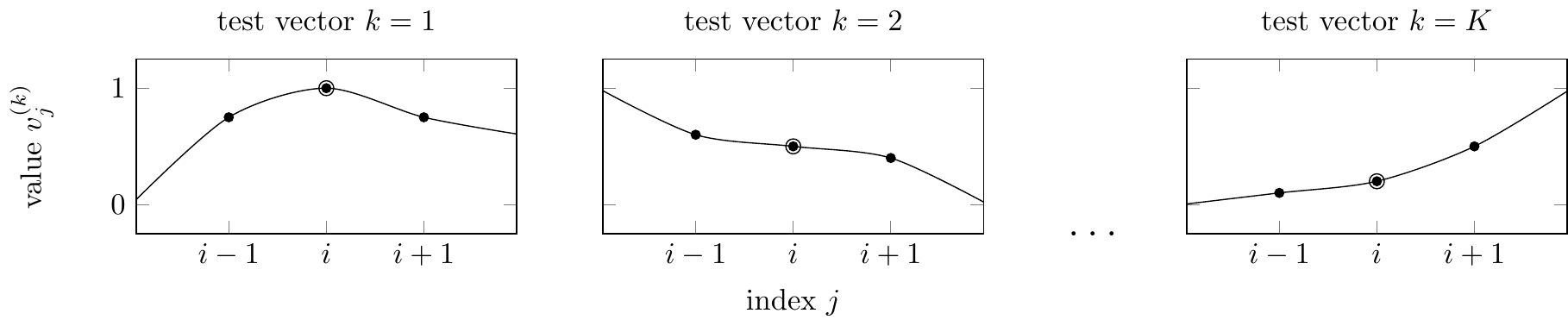}}
\[ \rightsquigarrow \text{ find } p_{i,i-1},p_{i,i+1} \text{ s.t. } \left(
\begin{bmatrix}
v^{(1)}_{\textcolor{black}{i}} \\
v^{(2)}_{\textcolor{black}{i}} \\
\vdots \\
v^{(K)}_{\textcolor{black}{i}}
\end{bmatrix}
-
p_{i,i-1}\begin{bmatrix}
v^{(1)}_{i-1} \\
v^{(2)}_{i-1} \\
\vdots \\
v^{(K)}_{i-1}
\end{bmatrix}
-
p_{i,i+1}\begin{bmatrix}
v^{(1)}_{i+1} \\
v^{(2)}_{i+1} \\
\vdots \\
v^{(K)}_{i+1}
\end{bmatrix}
\right)^2
\to \operatorname{min \, .}
\]
\centerline{\includegraphics[width=1\textwidth]{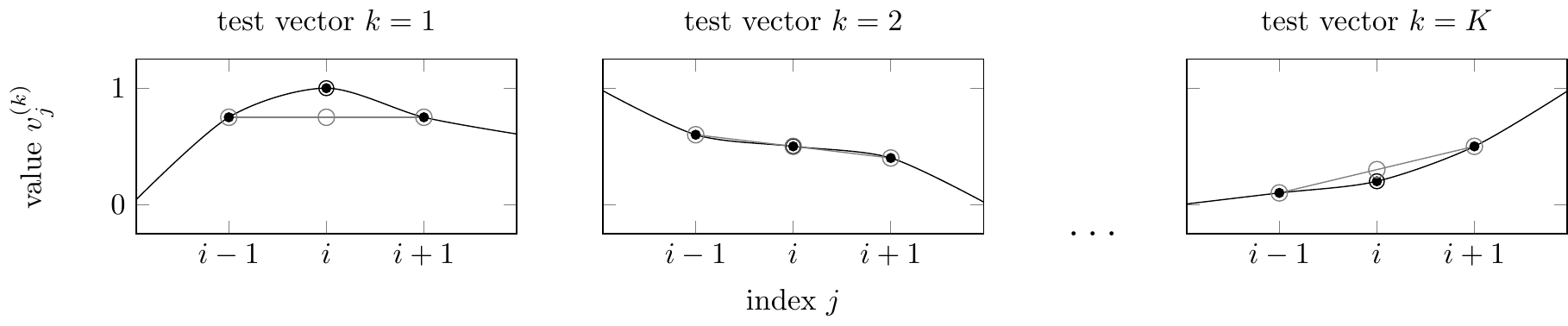}}
\caption{Example for \leastsquares{} interpolation from smooth error
\label{fig:lsinterpolation}}
\end{figure}

In order to explain the usefulness of \leastsquares{} interpolation in solving the coarsening problem we consider it in terms of a regression problem. Clearly, the interpolation weights $p_{ij}$ determined by \leastsquares{} interpolation can be thought of as an $\ell_2$ best regression fit to a set of observations of smooth error given by the entries of the test vectors. That is,~\eqref{eq:lsi} determines a best weighted $\ell_2$ fit to variable $i$ based on observations of smooth error made at selected (nearby) variables in $\mathcal{C}_{i}$ (cf.~\cref{fig:lsinterpolation}). While this regression fit is meaningful in case $K \gg |\mathcal{C}_{i}|$, i.e., we are in a data rich scenario, we run into a severe problem of overfitting in case $K \approx |\mathcal{C}_{i}|$. In the extreme case of $K \leq |\mathcal{C}_{i}|$ we obtain an exact fit of the (arbitrary) observations of smooth error which might lack generalizability; we refer to~\cite[Chapter 3]{HastTibsFrie2001} for a general introduction to regression in statistical learning. This is in accordance with earlier results reported for the bootstrap algebraic multigrid method where it was observed that a lack of data, i.e., test vectors, severly hampers the performance of the overall method~\cite{BranBranKahlLivs2011} or needs to be supplemented by implicit assumptions on the nature of algebraically smooth error~\cite{MantMcCoParkRuge2010}.

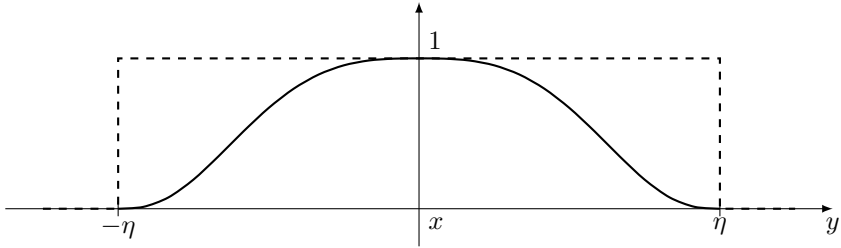
\begin{figure}
    \begin{center}
        \begin{tikzpicture}
            \draw[thin,-latex] (-5.5,0) -- (5.5,0)  node[below] {$y$};
            \draw[thin,-latex] (0,-.5) -- (0,2.75);
            \draw[thin] (-.1,2) -- node[above right] {$1$} (.1,2);
            \draw[thin] (-4,-.1) -- node[below] {$-\eta$} (-4,.1);
            \draw[thin] (4,-.1) -- node[below] {$\eta$} (4,.1);
            \node[below right] at (0,0) {$x$};
            \draw[thick,dashed] (-5,0) -- (-4,0) -- (-4,2) -- (4,2) -- (4,0) -- (5,0);
            \draw[domain=-4:4,smooth,variable=\x,thick]  plot ({\x},{2*(1-abs(\x/4)^3)^3});
        \end{tikzpicture}
    \end{center}
    \caption[]{Kernel functions in one dimension with $d(x,y) = \norm{x - y}$. Nearest-neighbor (\raisebox{.25em}{\resizebox{1em}{!}{\tikz{\draw[dashed] (0,0) -- (1,0);}}}) and tri-cube kernel (\raisebox{.25em}{\resizebox{1em}{!}{\tikz{\draw (0,0) -- (1,0);}}}).
    \label{fig:kernelfunctions}}
\end{figure}

Based on this interpretation of \leastsquares{} interpolation as a regression model of smooth error it should in principle be possible to extract information about an appropriate choice of $\mathcal{C}_{i}$ from the calculated regression coefficients, i.e., the interpolation weights $p_{ij}$. 
Due to the fact that with no a-priori information available on which variables might be important in the interpolation for variable $i$ we have to consider a \leastsquares{} fit from all other variables
\begin{equation}
    \sum_{k = 1}^{K} \omega_{k} \left(v_{i}^{(k)} - \sum_{j \neq i} p_{ij}v_{j}^{(k)}\right)^{2} \rightarrow \operatorname{min}.
\end{equation} This of course is ill-posed for the aforementioned reasons as long as one does not use an absurd amount of test vectors, i.e., $K > n$. The problem can be somewhat alleviated when assuming locality of interpolation and thus restricting the potential set of interpolation variables to some neighborhood\footnote{In accordance with an appropriate metric.} $\mathcal{N}_{i}$ of variable $i$, i.e., considering
\begin{equation}
    \sum_{k = 1}^{K} \omega_{k} \left(v_{i}^{(k)} - \sum_{j \in \mathcal{N}_{i}} p_{ij}v_{j}^{(k)}\right)^{2} \rightarrow \operatorname{min}.
\end{equation} A more general way of localizing the regression problem is to introduce what is known in statistical learning as a \textit{kernel} operator $K_{\eta}$ (cf.~the introduction to kernel regression in \cite[Chapter 6]{HastTibsFrie2001}) into the \leastsquares{} fit by,
\begin{equation}\label{eq:lsikernel}
    \sum_{k = 1}^{K} \omega_{k} \left(v_{i}^{(k)} - \sum_{j} p_{ij}K_{\eta}(i,j)v_{j}^{(k)}\right)^{2} \rightarrow \operatorname{min}.
\end{equation} The role of the kernel operator is to weigh the data depending on the distance of the data variable $j$ to the observation variable $i$. We consider two choices for $K_{\eta}$ in our tests, the nearest-neighbor kernel $K_{\eta}^{\rm nn}$ and the tri-cube kernel $K_{\eta}^{\rm tc}$, illustrated in~\cref{fig:kernelfunctions} and defined by:
\begin{equation}
\begin{array}{rcl}
    K^{\rm nn}_{\eta}(i,j) &=& \begin{cases} 1, & d(i,j) < \eta\\ 0, & \text{else} \end{cases} \\[.5em]
    K^{\rm tc}_{\eta}(i,j) &=& \begin{cases} (1 - |\frac{d(i,j)}{\eta}|^3)^3, & d(i,j) < \eta\\ 0, & \text{else} \end{cases}
    \end{array}\label{eq:kerneldef}
\end{equation} In these definitions $d$ denotes a distance function, e.g., graph distance or euclidean distance in case variable coordinates are known.

Yet, even with such a localization of the regression problem, the number of test vectors required would be prohibitively large as the number of potential interpolation variables $|\operatorname{supp}(K_{\eta}(i,j))| = |\mathcal{N}_{i}|$ significantly exceeds $|\mathcal{C}_{i}|,$ the number of variables to be used lateron in interpolation. In the context of algebraic distances or other related adaptive coarsening techniques this problem is circumvented by simply considering only binary relations, i.e., fits between any pairs of variables~\cite{BranBranKahlLivs2015b,LivnBran2012,NapoNota2016,Nota2010}, but one can easily imagine that it is possible to miss important group relations by only considering pairs. While early works of bootstrap algebraic multigrid introduce greedy strategies to choose $\mathcal{C}_{i}$ by adding one variable at a time with some success (cf.~\cite{BranBranKahlLivs2011}), but even a heuristic justification of this approach is questionable. A brute force approach, checking all sets of $m$ variables in the neighborhood to find the best set, is clearly too expensive as well.

\begin{figure}
    \begin{center}
        \begin{tikzpicture}
            \begin{scope}
                \draw[ultra thin, -latex] (-1.25,0) -- (2.25,0) node [right] {\small$q_{1}$};
                \draw[ultra thin, -latex] (0,-1.25) -- (0,2.75) node [above] {\small$q_{2}$};
                \draw[rotate=45,black!50] (2,0.8) ellipse (12pt and 6pt);
                \draw[rotate=45,black!50] (2,0.8) ellipse (24pt and 12pt);
                \draw[rotate=45,black!50] (2,0.8) ellipse (36pt and 18pt);
                \draw[rotate=45,black!50] (2,0.8) ellipse (48pt and 24pt);
                \draw[black] (0,0) circle (1);
                \filldraw (0.1736,0.9848) circle (2pt) node[below right] {$q$};
                \filldraw (0.85,1.98) circle (2pt) node[below right] {$\widehat{p}$};
            \end{scope}
            \begin{scope}[xshift={200}]
                \draw[ultra thin, -latex] (-1.25,0) -- (2.25,0) node [right] {\small$p_{1}$};;         
                \draw[ultra thin, -latex] (0,-1.25) -- (0,2.75) node [above] {\small$p_{2}$};;
                \draw[rotate=45,black!50] (2,0.8) ellipse (12pt and 6pt);
                \draw[rotate=45,black!50] (2,0.8) ellipse (24pt and 12pt);
                \draw[rotate=45,black!50] (2,0.8) ellipse (36pt and 18pt);
                \draw[rotate=45,black!50] (2,0.8) ellipse (48pt and 24pt);
                \draw[rotate=45,black] (-0.71,-0.71) rectangle (0.71,0.71);
                \filldraw (0,1) circle (2pt) node[above left] {$p$};
                \filldraw (0.85,1.98) circle (2pt) node[below right] {$\widehat{p}$};
            \end{scope}
        \end{tikzpicture}
    \end{center}
    \caption{The effect of penalization on a \leastsquares{} problem with solution $\widehat{p}$; (left) solution $q$ with $2$-norm penalization $\norm{q} < t$ (right) solution $p$ with $1$-norm penalization $\norm[1]{p} < t$. In the $2$-norm penalized solution we find $q_{1},q_{2}\neq 0$, but due to the shape of the $1$-norm unit cell we find $p_{2} = 0$ in the $1$-norm penalized solution.\label{fig:unitspheres}}
\end{figure}
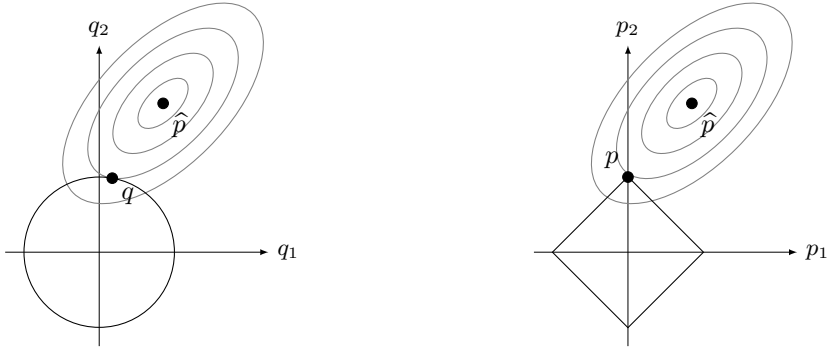

Thus we propose to extract the best set of variables by leveraging the similarity of the problem to a data regression problem. Especially under the premise that we do not want to increase the number of test vectors. One possible approach to sparsify $P$, which has been proposed in the machine learning context in~\cite{Tibs1996} is the introduction of an $\ell_{1}$ penalty term with corresponding penalty parameter $\lambda \in [0,\infty)$. This approach is known as least absolute shrinkage and selection operator (LASSO). To be more specific, the  \leastsquares{} problem~\eqref{eq:lsikernel} is changed to
\begin{equation}\label{eq:lasso}
    \sum_{k = 1}^{K} \omega_{k} \left(v_{i}^{(k)} - \sum_{j} p_{ij}K_{\eta}(i,j)v_{j}^{(k)}\right)^{2} + \lambda\norm[1]{p_{i}}\rightarrow \operatorname{min},
\end{equation} 
Thus, in essence $\lambda$ allows us to interpolate between the \leastsquares{} solution at $\lambda=0$ and $p_{i} = 0$ for $\lambda \rightarrow \infty$. Note, that this could also be achieved by penalizing the 2-norm $ \norm[1]{p_i}^2$, but penalizing the $1$-norm implicitely enforces sparsity of $p_{i}$. This observation can be motivated by stating \cref{eq:lasso} in the equivalent form
\begin{equation}\label{eq:lasso_t}
    \sum_{k = 1}^{K} \omega_{k} \left(v_{i}^{(k)} - \sum_{j} p_{ij}K_{\eta}(i,j)v_{j}^{(k)}\right)^{2} \rightarrow \operatorname{min} \quad \text{with} \quad \norm[1]{p_{i}} \leq t,
\end{equation} 
where $t$ being large corresponds $\lambda$ being small and vice versa. As illustrated in \cref{fig:unitspheres} the polyhedric shape of the $\ell_1$ unit sphere implicitely enforces zero entries in the $\ell_1$ penalized solution.

Now the idea is to filter out the most important interpolation variables by analyzing the regression weights $p_{ij}$ w.r.t.~$\lambda$ and in this way find an appropriate set of interpolatory variables for every grid point. Further postprocessing then ensures that a proper set of coarse variables is constructed based on this information (cf.~\cref{sec:larcoarsening}). In addition to~\eqref{eq:lasso} it is interesting to consider a sign-constrained version of this optimization problem
\begin{equation}\label{eq:lasso_signed}
    \sum_{k = 1}^{K} \omega_{k} \left(v_{i}^{(k)} - \sum_{j} p_{ij}K_{\eta}(i,j)v_{j}^{(k)}\right)^{2} + \lambda\norm[1]{p_{i}}\rightarrow \operatorname{min},\ s_jp_{ij} \geq 0
\end{equation} for pre-specified signs $s_{j} \in \{-1,1\}$.


For a specific $\lambda$, calculating $p_{ij}$ in~\cref{eq:lasso,eq:lasso_signed} is a non-linear optimization problem. Again this raises the question of feasibility of this approach. Fortunately, there is a way to solve the penalized \leastsquares{} problems~\cref{eq:lasso,eq:lasso_signed} for a particular set of $\lambda$ values without having to resort to non-linear optimization methods, which we introduce now.

\section{\Leastangle{} regression}\label{sec:lars} The method of \emph{\leastangle{} regression}, introduced in~\cite{EfroHastJohnTibs2004}, is a general approach for the solution of \leastsquares{} problems motivated by data regression. While it ultimately constructs the unrestricted \leastsquares{} solution its benefit in our situation stems from the fact that it can be shown that its intermediate stages solve the penalized problem at particular $\lambda$ values. The main idea of \LARS{} is to start with a zero initial guess for the coefficients and then ``switch on'' one \leastsquares{} coefficient at a time. 

To be specific, let us first introduce a generic \leastsquares{} problem to simplify notation:
\begin{equation}\label{eq:genleastsquares}
    \operatorname{min}_{x} \norm[2]{v - Wx},
\end{equation} where $v \in \mathbb{R}^{n}, x\in \mathbb{R}^{m}, W = \left[\begin{array}{c|c|c} w_{1} & \cdots & w_{m} \end{array}\right]\in \mathbb{R}^{n \times m}$ and we assume that the columns of $W$ are normalized.\footnote{Note, that when applying \LARS{} to the weighted \leastsquares{} fit described in~\cref{sec:lsinterpolation} normalization of the columns of $W$ is replaced by weighting with $\omega_{k}$. Further $n$ and $m$ correspond to the number $K$ of test vectors and the support of the kernel function $K_{\eta}$, respectively.} Given $x$, not necessarily the solution of~\cref{eq:genleastsquares}, the \textit{residual} of this \leastsquares{} problem is defined by $r = v - Wx$. Further we define the set $\mathcal{A} = \{i, x_{i} \neq 0\}$ of active variables and a notion of similarity of the columns $w_{i}$ and the current residual $r$,
\[
\rho_{i} := {|\innerprod{w_{i}}{r}|}
\] This coincides up to scaling by $\norm{r}$ with the cosine of the \emph{angle} between $r$ and columns $w_{i}$ of $W$ and can be interpreted as the absolute \emph{correlation} of these vectors.

In each iteration of \LARS{}, starting from a zero initial guess, i.e., $x=0$ and $\mathcal{A} = \emptyset$, \LARS{} chooses the variable $\ihat \notin \mathcal{A}$ with largest correlation $\rho_\ihat$, i.e., least angle, and adds it to the active set $\mathcal{A}$. Denoting by $d_{\mathcal{A}}$ the solution of the \leastsquares{} problem which is restricted to the current set of active variables $\mathcal{A}$, i.e.,
\begin{equation}\label{eq:activels}
    \operatorname{min}_{d_\mathcal{A}} \norm{r - W_{\cdot,\mathcal{A}}d_{\mathcal{A}}},
\end{equation} we find that $x_{\mathcal{A}}+d_{\mathcal{A}}$ solves the original \leastsquares{} problem~\cref{eq:genleastsquares} restricted to variables in $\mathcal{A}$. \CLARS{} now introduces a step size $\alpha \in (0,1]$ and updates the intermediate \leastsquares{} solution  by 
\begin{equation}
    \widetilde{x} = \begin{cases} x_{i}+\alpha \cdot d_{i} & i \in \mathcal{A}\\x_{i}&\text{else}\end{cases}
\end{equation} and correspondingly $\widetilde{r} = r - \alpha W_{\cdot,\mathcal{A}}d_{\mathcal{A}}$. Defining updated correlations by 
\[
\widetilde{\rho} = W^{H}\widetilde{r} = W^{H}r - \alpha \underbrace{W^{H}W_{\cdot,\mathcal{A}}d_{\mathcal{A}}}_{=:\ \mu} = \rho - \alpha \cdot \mu
\] the step size $\alpha$ in \LARS{} is then chosen as the smallest positive $\alpha$ fulfilling any of the following conditions
\begin{enumerate}[(i)]
    \item $\widetilde{x}_{i} = 0,\ i \in \mathcal{A}$,
    \item $|\widetilde{\rho}_{i}| = |\widetilde{\rho}_{j}|,\ i \in \mathcal{A},\ j\notin \mathcal{A}$
\end{enumerate}
In case $\alpha$ is chosen due to $\widehat{x}_{i} = 0$ for some $\icheck \in \mathcal{A}$, the $\icheck$-th variable is removed from the active set, i.e., $\mathcal{A} \rightarrow \mathcal{A}\setminus \{\ \icheck\ \}$, and a new solution of~\cref{eq:activels} is calculated. Dropping variables from the active set when they become zero ensures that the obtained solution is equivalent to a solution of the penalized \leastsquares{} problem~\cref{eq:lasso} for some $\lambda \in [0,\infty)$; cf.~\cite{EfroHastJohnTibs2004}. The corresponding $\lambda$ value of~\cref{eq:lasso} is not known in \LARS{}, but subsequent iterations correspond to decreasing $\lambda$ values.

If, on the other hand, $\alpha$ is chosen due to a condition of type (ii), the corresponding variable $\ihat \notin \mathcal{A}$ is added to $\mathcal{A}$ for the next iteration. This guarantees that all variables in the active set are tied in correlation with the residual at all times; cf.~\cite{EfroHastJohnTibs2004}. 

\begin{algorithm2e}[ht]
    \SetAlgorithmStyle
    \caption{\Leastangle{} regression}
    \label{alg:leastangleregression}
    \KwData{$v \in \mathbb{R}^{n}, W = \left[\begin{array}{c|c|c}w_{1} & \cdots & w_{m}\end{array}\right] \in \mathbb{R}^{n\times m}$}\medskip
    Initialize $x = 0, \mathcal{A} = \emptyset$ and $\mathtt{rewind} = \mathtt{false}$ \;
    \While{$|\mathcal{A}| \leq \operatorname{min}(n,m)$}{
        $r = v - Wx$\;
        $\rho = W^{H}r$\;
        $\ihat = \operatorname{argmax}_{j \notin \mathcal{A}} | \rho_{j} | $\label{alg:defcorrelation}\;
        \uIf{$\mathtt{rewind} = \mathtt{false}$}{
            $\mathcal{A} = \mathcal{A} \cup \{\ \ihat\ \}$ \label{alg:activeset}\;
        }\Else{
            $\mathtt{rewind} = \mathtt{false}$\;
        }
        $d_{\mathcal{A}} = \operatorname{argmin}_{z_{\mathcal{A}}} \norm{r - W_{\cdot,\mathcal{A}}z_{\mathcal{A}}}$ \label{alg:defd} \;
        \For{$i \in \mathcal{A}, i \neq \ihat$}{
            $\alpha_{i} = -\frac{x_{i}}{d_{i}}$\;
        }
        $\widecheck{\alpha} = \operatorname{min}_{i\in \mathcal{A}}\{\alpha_{i} \geq 0\}$,\ $\icheck = \operatorname{argmin}_{i\in \mathcal{A}}\{\alpha_{i} \geq 0\}$\;
        $\mu = W^{H}W_{\cdot,\mathcal{A}}d_{\mathcal{A}}$\;
        \For{$j \notin \mathcal{A}$}{
            $(\alpha_{j}^{\prime},\alpha_{j}^{\prime\prime}) = \left(\frac{\rho_{i}-\rho_{j}}{\mu_{i}-\mu_{j}}, \frac{\rho_{i}+\rho_{j}}{\mu_{i}+\mu_{j}}\right),\ i \in \mathcal{A}$ \label{alg:intersections}\;
        }
        $\widehat{\alpha} = \operatorname{min}_{j\notin \mathcal{A}}\{\alpha_{j}^{\prime}, \alpha_{j}^{\prime\prime} \geq 0\}$\;
        \uIf{$\widecheck{\alpha} < \widehat{\alpha}$}{ \label{alg:alphadesc}
            $\mathtt{rewind} = \mathtt{true}$\;
          $x_{\mathcal{A}} = x_{\mathcal{A}} + \operatorname{min}\left\{1,\widecheck{\alpha}\right\} \cdot d_{\mathcal{A}}$\;
        }\Else{
            $x_{\mathcal{A}} = x_{\mathcal{A}} + \operatorname{min}\left\{1,\widehat{\alpha}\right\} \cdot d_{\mathcal{A}}$\;
        }
        \If{$\mathtt{rewind} = \mathtt{true}$}{
            $\mathcal{A} = \mathcal{A}\setminus \{\ \icheck\ \}$\;   
        }
    }
\end{algorithm2e}
The \LARS{} algorithm is summarized in~\cref{alg:leastangleregression}. In situations, where $n \ll m$ its computational complexity scales only linearly in the large dimension. This is of particular importance, when we apply the algorithm to~\cref{eq:lasso} where we find for the number of test vectors $K$ and the support of the kernel function $K_{\eta}$ that $K \ll \operatorname{supp}(K_{\eta})$. That is, \LARS{} scales only linearly with the size of the considered neighborhood.

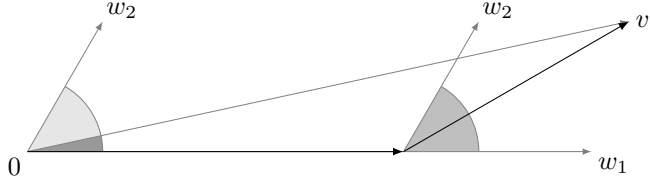
\begin{figure}
    \begin{center}
        \begin{tikzpicture}
            \fill[fill=black!40!white] (0,0) -- (1,0) arc (0:12.2:10mm) -- cycle;
            \draw[very thin,black!50!white] (1,0) arc (0:12.2:10mm);
            
            \fill[rotate=12.2,fill=black!10!white] (0,0) -- (1,0) arc (0:47.8:10mm) -- cycle;
            \draw[very thin,rotate=12.2,,black!50!white] (1,0) arc (0:47.8:10mm);
            
            \fill[fill=black!25!white] (5,0) -- (6,0) arc (0:60:10mm) -- cycle;
            \draw[very thin,black!50!white] (6,0) arc (0:60:10mm);
            
            \draw[-latex,very thin,black!50!white] (0,0) node[below left=-1pt,black] {$0$} -- (1,{1*sqrt(3)}) node[above right=-2pt,black] {$w_{2}$};
            \draw[-latex,very thin,black!50!white] (0,0) -- (7.5,0) node[below right=-1pt,black] {$w_{1}$};
            \draw[-latex,name path=initial target,very thin,black!50!white] (0,0) -- (8,{sqrt(3)}) node[right=-1pt,black] {${v}$};
            \draw[-latex,very thin,black!50!white] (5,0) -- (6,{1*sqrt(3)}) node[above right=-2pt,black] {$w_{2}$};
            
            \draw[-latex] (0,0) -- (5,0);
            \draw[-latex] (5,0) -- (8,{sqrt(3)});
        \end{tikzpicture}
    \end{center}
    \caption[Solution of a two-dimensional \leastsquares{} problem by \LARS{}]
    {Solution of a two-dimensional \leastsquares{} problem by \LARS{}. In the first iteration \LARS{} proceeds along $w_{1}$ until the residual bisects the angle between $w_{1}$ and $w_{2}$. In the second and final iteration \LARS{} proceeds along the bisector.\label{fig:2dlar}\label{fig:largeometric}}
\end{figure}

\begin{figure}[htb]
    \includegraphics[width=1\textwidth]{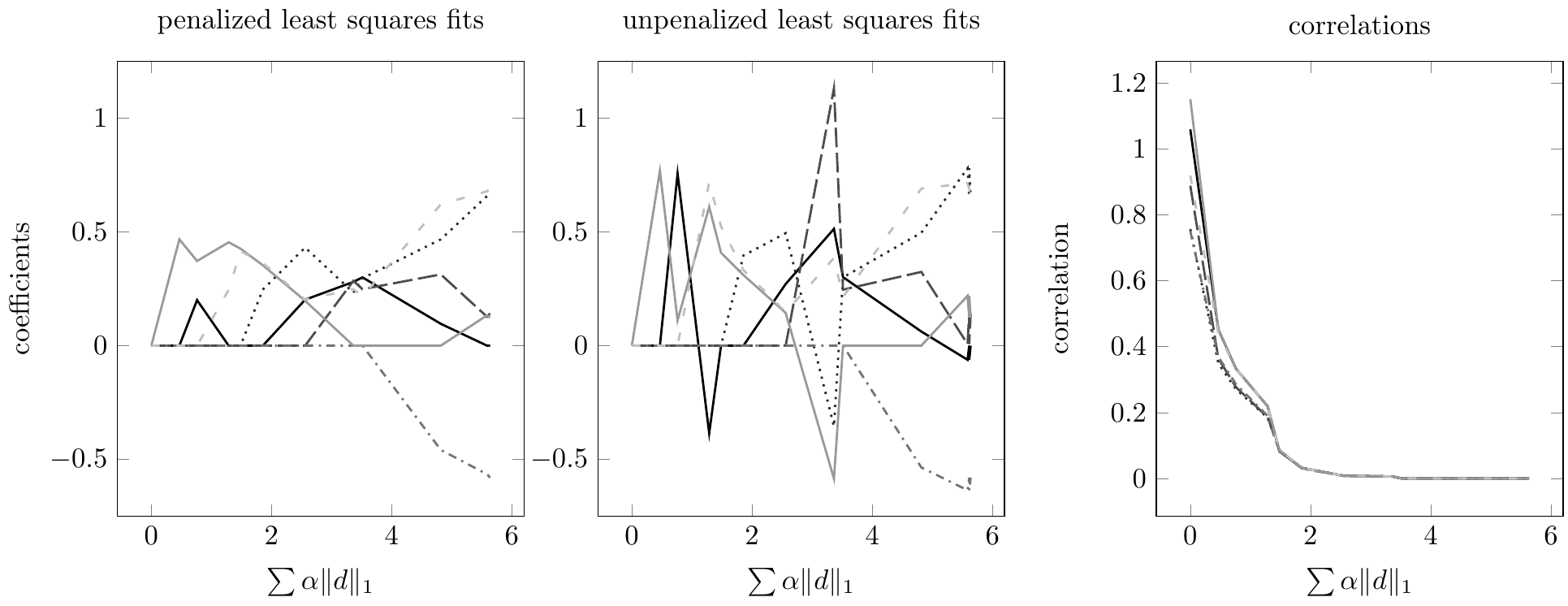}
    \caption{\label{fig:larexample} An example of the coefficient traces of a \LARS{} iteration with respect to accumulated $\ell_1$ coefficient changes along the \LARS{} iteration; 
    (left) penalized \leastsquares{} coefficients (\LARS{} $\alpha$); (middle) unpenalized \leastsquares{} coefficients ($x + d_{\mathcal{A}}$ in \LARS); (right) correlations $\rho$.}
\end{figure}

In~\cref{fig:largeometric} we illustrate the \LARS{} iteration in the case of $W \in \linspace[n\times 2]{R}$. The property of tied correlations transfers geometrically to choosing a step size such that the remaining line to the projection $\widehat{v}$ of the measurement $v$ bisects the angle between the currently active direction $w_{1}$ and the inactive $w_{2}$. 
In~\cref{fig:larexample} \LARS{} coefficients and the respective correlations are shown for a simple example with $n=8$ and $m=6$. That is, we calculate up to $6$ \LARS{} coefficients based on $8$ vectors. Note, that we have chosen $n\geq m$ in this case to illustrate a complete \LARS{} trajectory that ends with a well-defined \leastsquares{} solution.


\paragraph{\Leastangle{} regression with sign constraint}
With a small modification \LARS{} is able to solve a sign-constrained version of the penalized \leastsquares{} problem~\cref{eq:lasso_signed}. In this the signs of the calculated coefficients are tied to the signs of their initial correlations.
In order to ensure these tied signs we simply modify the data matrix $W$ by scaling it with the signs of the initial correlations $\rho^{(0)} = W^{H}v$, i.e.,
\begin{equation}\label{eq:monotonicLARSI}
    \widehat{W} = W\operatorname{diag}(\operatorname{sign}(\rho^{(0)})).
\end{equation} This guarantees that the modified initial correlations $\widehat{W}^{H}v$ are now all positive. As the coefficients entering the active set in~\cref{alg:activeset} of~\cref{alg:leastangleregression} over the course of \LARS{} have the same sign as their correlation (cf.~\cite{EfroHastJohnTibs2004}) we simply have to ensure that only coefficients with positive correlations enter the active set. For this we need two additional modifications of~\cref{alg:leastangleregression}. First, we have to omit considering the absolute values in~\cref{alg:defcorrelation}, i.e., replacing the statement by
\[
\ihat = \operatorname{argmax}_{j \notin \mathcal{A}} \rho_{j} \, ,
\] to ensure that the chosen maximal correlation is positive. This can be guaranteed if the \LARS{} iteration is stopped once all correlations reached or crossed zero. 
Second, instead of considering
\[
|\widetilde{\rho}_{j}| = |\widetilde{\rho}_{i}|,\ j \notin \mathcal{A}, i \in \mathcal{A}
\]
in~\cref{alg:intersections} we only consider $\widetilde{\rho}_{j} = \widetilde{\rho}_{i}$ to determine $\widehat{\alpha}$. That is, we replace~\cref{alg:intersections} by
\begin{equation}\label{eq:monotonicLARSIII}
    {\alpha}_{j} = \frac{\rho_{i}-\rho_{j}}{\mu_{i}-\mu_{j}} \, .
\end{equation} Once the algorithm terminates we obtain the coefficients with the correct sign by reverting the sign change on the data matrix, i.e., calculating
\begin{equation}
    x \leftarrow \operatorname{diag}(\operatorname{sign}(\rho^{(0)})) x \, .
\end{equation}

\paragraph{Stopping criteria}
In order to not only provide an unpenalized \leastsquares{} fit, which it does once all variables are active, \LARS{} requires a suitable stopping criterion. In the context of the coarsening problem two choices come to mind. 

First, due to the fact that the data in the coarsening problem should be highly correlated as representations of algebraically smooth error, one can use the correlation of the inactive variables to decide when to stop the iteration. As illustrated in~\cref{fig:larexample} the correlations quickly become very small and a stopping criterion based on the largest inactive correlation seems appropriate. In addition, a stopping criterion in terms of remaining correlation can be interpreted in terms of the amount of information absorbed in the model: If the remaining observations are already well approximated by the model, there is no further need to fit them. 

Second, specifically for the coarsening problem, one can stop \LARS{} by using either the cardinality of $\mathcal{A}$, i.e., the cardinality of $\mathcal{C}_{i}$ in \leastsquares{} interpolation. Given the fact that the number of elements of $\mathcal{A}$ is non-monotonic one cannot stop the iteration exactly at a prescribed cardinality, but rather back-track the last occurrence of a certain cardinality after stopping at a high enough cardinality to avoid missing the last occurrence. The largest number of \LARS{} iterations is naturally bounded by either $n$ or $m$, but oftentimes stopping the iteration with $|\mathcal{A}|$ of twice the prescribed cardinality almost always allows the extraction of the corrected $\mathcal{C}_{i}$ set even if the natural limit of iterations is not reached. Any \LARS{} solution with a cardinality close to the number of available test vectors should be treated with utmost care due to the potential problem of overfitting. The left and the center panel of \cref{fig:larexample} illustrate this problem as the \leastsquares{} solution contains negative weights in order to balance the sum of weights. Such fits are very sensitive under small changes in the input data.


\section{\Leastangle{} regression coarsening}\label{sec:larcoarsening}
Having introduced the penalized \leastsquares{} problem and the \LARS{} method to solve it, we now construct our adaptive coarsening algorithm. Our approach is to use the coefficients of \LARS{} to define a notion of strength of connection. 


The first step in this process is to determine penalized and unpenalized \leastsquares{} coefficients $p_{ij}$ by \LARS{} for all variables $i\in \Omega$ for a given diameter $\eta$ of the kernel function and choice of distance function $d$. In addition, one has to specify a stopping criterion for \LARS{} as discussed in the previous section. 
After truncation of coefficients below a specified threshold, the remaining coefficients define a directed graph of strong connections. As a first guess at an appropriate set of coarse variables, an independent set of this graph is computed (cf.~\cref{alg:indset}), where the variables are chosen according to the importance measure
\begin{equation}
    \sigma_j = \sum_{i, j \in \mathcal{N}_{i}} |p_{ij}| \, .
\end{equation} That is, the priority of variables is determined by the weight of strong couplings contributed to other variables.

\begin{algorithm2e}[ht]
    \SetAlgorithmStyle
    \caption{Independent set coarsening with ordering by importance.\label{alg:indset}}
    \KwData{Sets $C_{i}$ and \leastsquares{} coefficients $p_{ij}$ for all $i$, threshold $\theta$}
    \KwResult{Coarse variable set $\mathcal{C}$}
    \For{$i = 1,\ldots,n$}{
        \For{$j \in \mathcal{C}_{i}$}{
            \If{$\frac{|p_{ij}|}{\operatorname{max}_{j}|p_{ij}|} < \theta$}{
                $p_{ij} = 0$\;
                $\mathcal{C}_{i} = \mathcal{C}_{i}\setminus \{j\}$\;
            }
        }
    }
    Initialize importance scores $\sigma_{i} = 0,\ i = 1,\ldots,n$\;
    \For{$i = 1,\ldots,n$}{
        \For{$j \in \mathcal{C}_{i}$}{
            $\sigma_{j} = \sigma_{j} + |p_{ij}|$\;
        }
    }
    Initialize $\mathcal{C} = \emptyset,\ \mathcal{B} = \{1,\ldots,n\}$\;
    \While{$\mathcal{B} \neq \emptyset$}{
        $i^{\star} = \operatorname{argmax}_{i} \{\sigma_{i}, i \in \mathcal{B}\}$\;
        $\mathcal{C} = \mathcal{C} \cup \{i^{\star}\}$\;
        $\mathcal{B} = \mathcal{B} \setminus \{i^{\star}\}$\;
        \For{$i \in \mathcal{B}, p_{ii^{\star}} \neq 0$}{
            $\mathcal{B} = \mathcal{B} \setminus \{i^{\star}\}$\;
        }
    }
\end{algorithm2e}

Once an initial set of coarse variables is known, another pass of \LARS{} for all variables $i \in \Omega \setminus \mathcal{C}$ is carried out, but now the kernel function is restricted to variables in $\mathcal{C}$. This determines a first set of interpolatory variables $\mathcal{C}_{i}$ and a first set of \leastsquares{} interpolation weights $p_{ij}$.

Due to the fact that the independent set might yield inconsistent sets of coarse variables, which in turn might lead to inefficient interpolation we further modify $\mathcal{C}$ by alternating between the calculation of interpolatory sets $\mathcal{C}_{i}$ and interpolation weights $p_{ij}$ in the aforementioned way and a \emph{maximal volume} correction. 

The idea of the maximal volume correction is to find a choice of coarse variables such that no interpolation weight is larger in absolute magnitude than $1$. This is achieved by subsequently finding the largest entry $|p_{ij}|$ of interpolation and  swapping variable $j \in \mathcal{C}$ with variable $i \in \mathcal{F}$ in case $|p_{ij}| > 1$. This process has been introduced in~\cite{GoreOselSavoTyrtZama2010,Knut1985} to calculate well conditioned bases of linear spaces and maximal volume submatrices, and it can be shown that this process converges, requires only rank-$1$ updates of the interpolation weight matrix and yields the desired result, i.e.,  interpolation weights of modulus smaller one. As the maximal volume correction changes the sparsity of interpolation, we afterwards rerun \LARS{} to calculate new sets $\mathcal{C}_{i}$ and interpolation weights $p_{ij}$ for the variables in $\mathcal{F}$ with modified sparsity. The whole process is repeated until no corrections are done in the maximal volume part or a pre-set number of iterations has been reached.

\begin{algorithm2e}[ht]
    \SetAlgorithmStyle
    \caption{Maximal volume correction\label{alg:maxvol}}
    \KwData{interpolation $P$, coarse variable set $\mathcal{C}$}
    \KwResult{coarse variable set $\mathcal{C}$}
    \Repeat{$\max_{i,j} |p_{ij}| \leq 1$}{
        find $k \notin \mathcal{C}, \ell \in \mathcal{C}$ such that $|p_{k\ell}| = \max_{i,j}|p_{ij}| > 1$\;
        add $k$ to $\mathcal{C}$, remove $\ell$ from $\mathcal{C}$\;
        update entries of $P$\;
    }
\end{algorithm2e}

In case that after a \LARS{} call a variable $j\in \mathcal{C}$ turns out to not interpolate to any variable $i \in \mathcal{F}$, $j$ is added to the $\mathcal{F}$ variables and \LARS{} is used to determine $\mathcal{C}_{i}$ and $p_{i}$ for these new $\mathcal{F}$ variables.

\begin{algorithm2e}[ht]
    \SetAlgorithmStyle
    \caption{\Leastangle{} regression coarsening\label{alg:larcoarsening}}
    \KwData{test vectors $v^{(1)},\ldots,v^{(K)}$}
    \KwResult{coarse variable set $\mathcal{C}$ and \leastsquares{} interpolation $P$}
    \For{all variables $i \in \Omega$}{
        calculate $\mathcal{C}_{i}, p_{i}$ by \texttt{\LARS{}}\;
        apply thresholding to $p_{i}$ and change $\mathcal{C}_{i}$ accordingly\;
    }
    Calculate $\mathcal{C}$ by \texttt{independent set}\;
    \Repeat{no maximal volume correction occurred}{
        calculate $\mathcal{C}_{i}, p_{i}$ by \texttt{\LARS{}} with $K_{\eta}\cap\mathcal{C}$\;
        \For{any variable $j \in \mathcal{C}$ with $j \notin \mathcal{C}_{i}, i \notin \mathcal{C}$}{
            remove $j$ from $\mathcal{C}$\;
            calculate $\mathcal{C}_{i}, p_{i}$ by \texttt{\LARS{}} with $K_{\eta}\cap\mathcal{C}$\;
        }
        obtain updated $\mathcal{C}$ set by \texttt{maximal volume correction}\;
    }
\end{algorithm2e}

\section{Numerical experiments}\label{sec:numericalresults}


We now show some tests of our MATLAB implementation of the \LARS{} coarsening approach that we ran on finite element discretizations of Poisson's equation with and without anisotropy. We thus consider the partial differential equation 
\begin{equation} \label{eq:poisson_univ}  - \left( c_1 \frac{\partial^2}{\partial x^2} + c_2 \frac{\partial^2}{\partial y^2} + c_3 \left( \frac{\partial}{\partial x} \frac{\partial}{\partial y} + \frac{\partial}{\partial y}  \frac{\partial}{\partial x} \right) \right) u = f \, ,\end{equation}
where $c_1 = c_2 = 1$ and $c_3 = 0$ yields Poisson's equation
\begin{equation} \label{eq:poisson}
- \left( \frac{\partial^2}{\partial x^2} + \frac{\partial^2}{\partial y^2} \right) u = f \, ,
\end{equation}
and for
\begin{equation*}
c_1 = \cos(\alpha)^2 + \varepsilon\cdot\sin(\alpha)^2, \quad c_2 = \sin(\alpha)^2 + \varepsilon\cdot\cos(\alpha)^2
\end{equation*}
and
\begin{equation*}
c_3 = \frac{(1-\varepsilon)}{2}\cdot\sin(2\cdot\alpha)
\end{equation*}
we obtain Poisson's equation with anisotropy in direction $\alpha$. The underlying domain $\Omega = \{ x \in \mathbb{R}^2 \, : \, \|x\| \leq 1 \}$ is the unit disc and the finite element discretization (that remains unchanged for all tests in this section) has been obtained from MATLAB's \texttt{pdetool} using linear elements and Dirichlet boundary conditions $u |_{\partial\Omega} = 0$.



\subsection{\CLARS{}}

\begin{figure}[tbh]
\centerline{\includegraphics[width=1\textwidth]{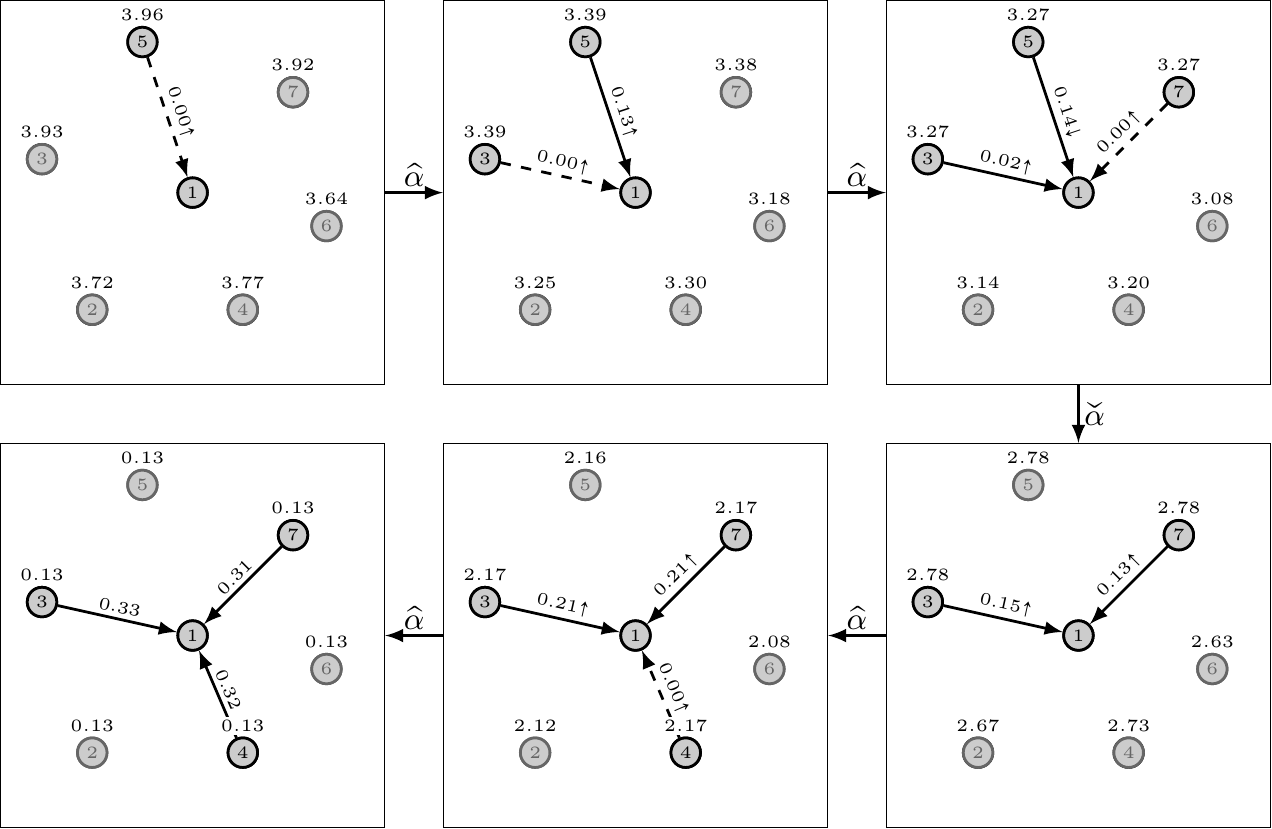}}
\caption{Illustration of \LARS{} with a kernel range of one edge. Each panel depicts the state in \cref{alg:alphadesc} of \cref{alg:leastangleregression} in consecutive iterations, i.e., right before updating the (penalized) coefficients $x$. For each node we report $|\rho_i|$. We highlight the current active set $\mathcal{A}$. Edge values denote the current penalized \LARS{} coefficients $x_{i}$ and the small arrow indicates the sign of $d_{i}$, i.e., the direction in which the coefficient is about to be changed. The connection between panels is marked either by $\widehat{\alpha}$ or $\widecheck{\alpha}$ depending on which of the two is smaller in the condition in~\cref{alg:alphadesc}.
\label{fig:lars_local}}
\end{figure}

To start our numerical tests, we would like to highlight what makes \LARS{} advantageous in the coarsening process of algebraic multigrid. To do so we consider Poisson's equation discretized with $h_{\rm max}=0.1$. This yields a matrix $A\in \mathbb{R}^{480 \times 480}$. We use $K=16$ test vectors, each of them smoothed with $4$ steps of Gauss-Seidel. For convenience we use the notation of \cref{alg:leastangleregression}.

In the first \LARS{} iteration variable $5$ is added to $\mathcal{A}$ as it has the largest correlation with variable $1$. The weight $x_5$ is increased and variable $3$ is added next as the correlations of variables $3$ and $5$ are now equal. In the third iteration variable $7$ is added. While $x_7$ increases, $x_5$ decreases with the next update. The step size $\alpha$ in this iteration is chosen such that variable $5$ is dropped $(\widecheck{\alpha} < \widehat{\alpha})$. That is, variable $5$ is obsolete in case variables $3$ and $7$ are in $\mathcal{A}$ which intuitively makes sense. Afterwards variable number $4$ is added, resulting in a caliber three interpolation set which is a geometrically balanced choice and the last occurence of a caliber three set in the \LARS{} iterations for this variable.

In this way~\cref{fig:lars_local} illustrates why we deem \LARS{} to be advantegous compared to other adaptive coarsening approaches. Binary approaches that only consider the relation between pairs of variables would not be able to devalue variable $5$ in presence of variables $3$ and $7$, but would treat all three connections as equal. On the other hand, it describes an almost greedy way, guided by the $\ell_{1}$ constraint \leastsquares{} problem, to construct interpolation sets $\mathcal{C}_{i}$ with more than one interpolation point. It is able to do so without testing all possible combinations and offers the possibility to adjust key parameters of the coarsening such as the caliber of interpolation in an adaptive fashion based on quantities such as the correlation.

\subsection{Parameter study} We first review the available parameters and group them into \emph{kernel}, \emph{\LARS} and \emph{coarsening} parameters.

\begin{itemize}
\item \textbf{Kernel parameters} include the \emph{distance measure} $d$, the \emph{kernel function} and the \emph{kernel radius} $\eta$ in~\cref{eq:kerneldef}.
\begin{itemize}
\item \emph{distance measure:} Even though information about the coordinates of each variable is available, we choose to use the graph distance.
\item \emph{kernel function:} We use both nearest-neighbor and tri-cube kernel and specify this for each test individually.
\item \emph{kernel radius:} For all tests we used $\eta = 4$.
\end{itemize}
\item \textbf{\CLARS\  parameters} consist of a \emph{correlation threshold}, a \emph{relative strength threshold}, a \emph{caliber threshold} and a flag for the use of the \emph{sign constraint}.
\begin{itemize}
\item \emph{correlation threshold:} Stopping criterion based on the current correlation of inactive variables in~\cref{alg:leastangleregression}.
\item \emph{relative strength threshold:} Truncation of penalized/unpenalized regression coefficients smaller than threshold times coefficient of largest absolute value.
\item \emph{caliber threshold:} Chooses penalized/unpenalized regression coefficients in~\cref{alg:leastangleregression} of the iteration with last occurrence of $|\mathcal{A}|$ equal to threshold.
\item \emph{sign constraint:} Switches between the lasso \leastsquares{} problem~\cref{eq:lasso} and its sign constrained version~\cref{eq:lasso_signed}.
\end{itemize}
\item \textbf{Coarsening parameter} specify the \emph{number of maximal volume and \LARS{} loops}.
\end{itemize}

In order to study the behavior of the \LARS{} coarsening approach with respect to these parameters and to come up with well-founded default choices we first study their individual influence on a small sample problem. To this end we use \texttt{pdetool}
with $h_{\rm max} = 0.1$ to generate a matrix $A\in \mathbb{R}^{480 \times 480}$ for the Poisson equation~\cref{eq:poisson}. 

We construct $2$- and $3$-grid methods using $K = 16$ initially random i.i.d.~$N(0,1)$-distributed test vectors and apply $4$ Gauss-Seidel iterations to each to obtain a set of smoothed test vectors. Reported convergence results correspond to a V$(1,1)$-cycle with symmetric\footnote{Forward sweep on the fine-to-coarse and backward sweep on the coarse-to-fine pass.} Gauss-Seidel smoothing. We use the same parameter sets for both coarsening steps and apply $4$ smoothing iterations to the restricted test vectors on the intermediate grid in the construction of the 3-grid method. All parameters that are not explicitly set in the following tests are chosen such that they do not influence the method.

\begin{figure}[tbh]
\centerline{\includegraphics[width=1\textwidth]{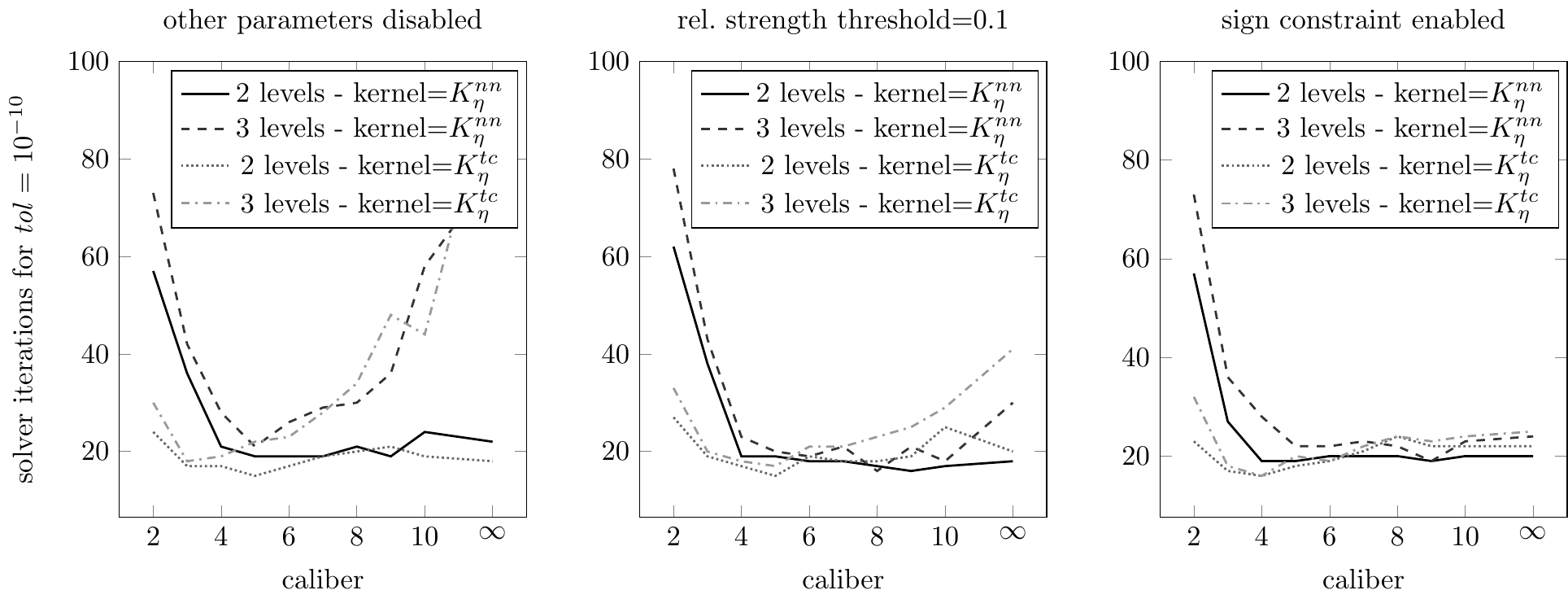}}
\caption{Study of iteration count as a function of the caliber threshold. 
\label{fig:caliber_study}}
\end{figure}
The first parameter we study is the influence of the \emph{caliber threshold} on the convergence of the $2$- and $3$-grid method. As illustrated in~\cref{fig:caliber_study} in the left plot we see on one hand that there exists a lower bound on the caliber threshold that is required to obtain a rapidly converging method. On the other hand the method becomes unstable, i.e., its convergence degradates when going from a two to a three grid method, for large calibers which might be explained by the problem of overfitting the data as mentioned in~\cref{sec:lsinterpolation}. This problem is especially apparent in the $3$-grid method where no stable plateau is visible. Using a relative strength threshold or the sign constraint version as depicted in the middle and right plot, respectively, cures this instability at large caliber almost completely and a stable plateau arises. 
\begin{figure}[bht]
\centerline{\includegraphics[width=1\textwidth]{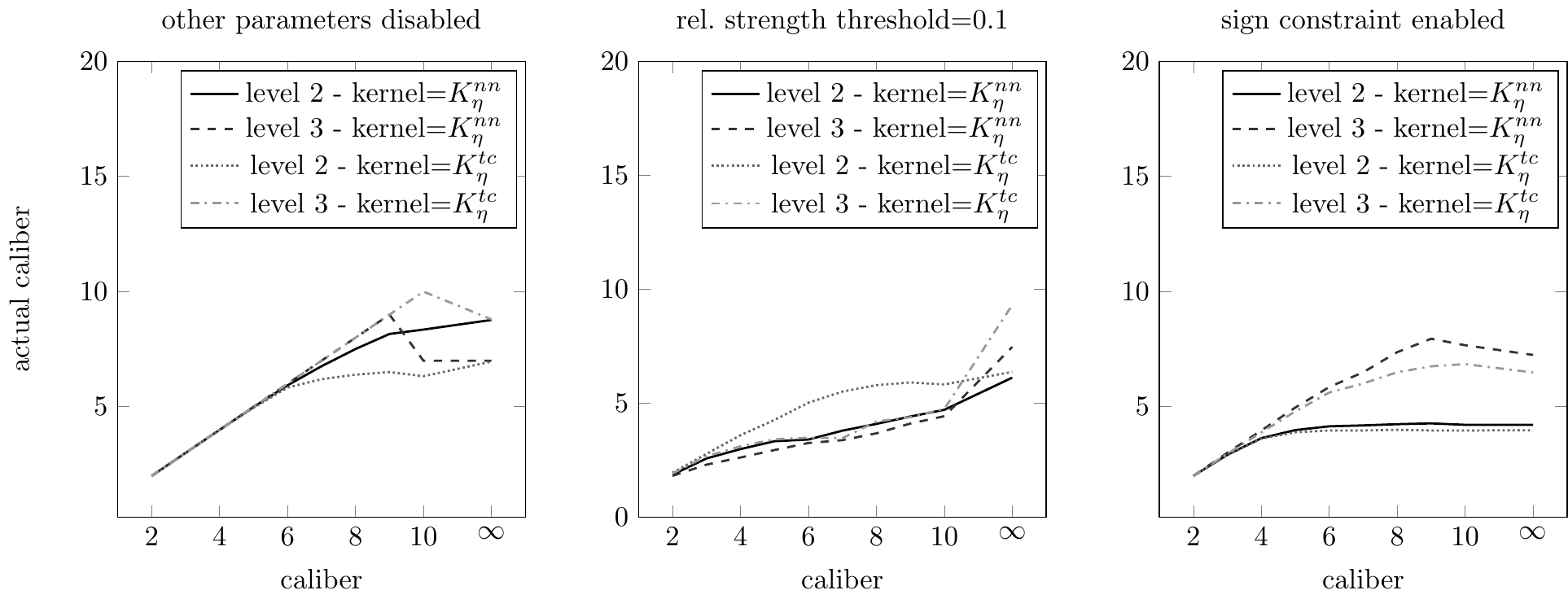}}
\caption{Study of actual caliber as a function of the caliber threshold.
\label{fig:act_caliber_study}}
\end{figure}
This is unsurprising when we combine the observed behavior with the plots in~\cref{fig:act_caliber_study}, where we plot the actual caliber for the same tests. The actual calibers for the $3$-grid cases denote the calibers for the coarsening on the intermediate grids that are identical to the coarse grids of the respective $2$-grid method. While the actual caliber grows almost linearly in the left plot, where no additional stopping criterion is used, we can see that both additional parameters effectively stabilize the actual caliber in a range small enough to not cause overfitting and thus reduce the instability of the method. Based on these findings, we set the default choice for the caliber threshold to $3$. With respect to the choice of kernel function no dramatic difference between nearest-neighbor and tri-cube is visible, with a minor advantage of tri-cube in most cases, which again is not surprising as it enforces stricter locality of the coarsening.

\begin{figure}[tbh]
\centerline{\includegraphics[width=1\textwidth]{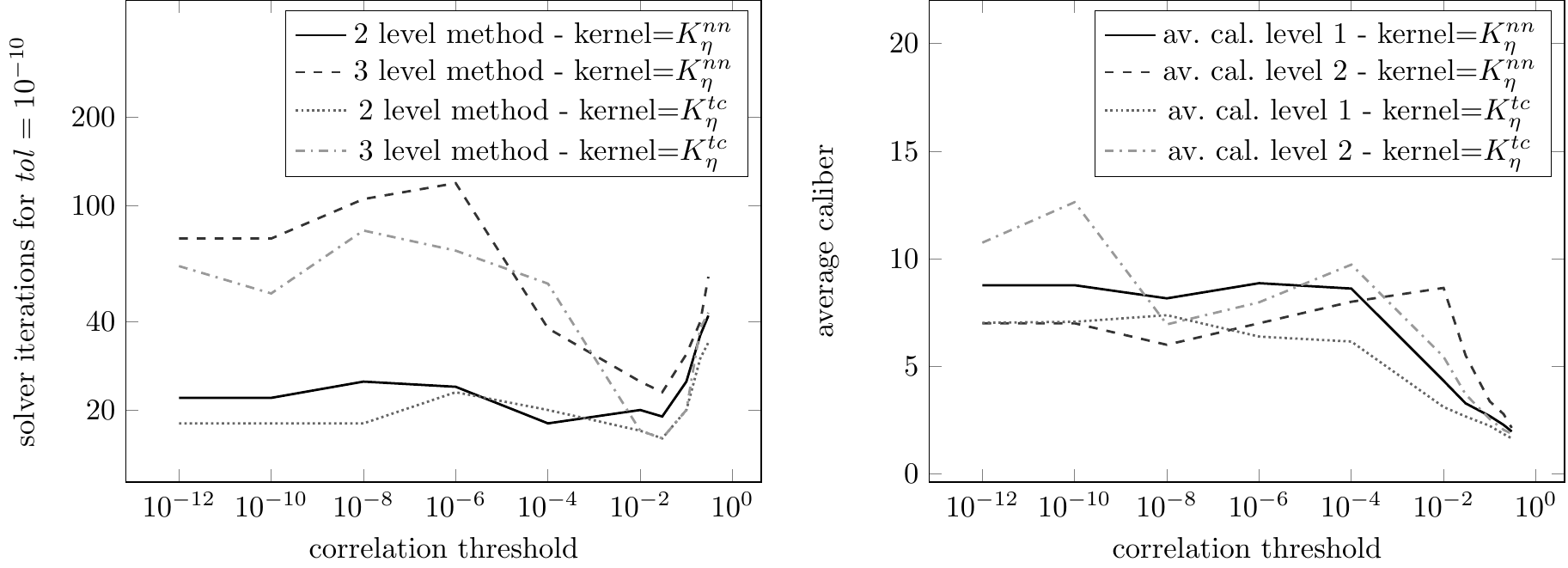}}
\caption{Study of iteration count (left) and actual caliber (right) with respect to the correlation threshold.\label{fig:mincorr_study}}
\end{figure}
Next, we consider the influence of the correlation threshold on convergence in~\cref{fig:mincorr_study}. As one can see, the convergence behavior of the $2$-grid method is very robust with respect to this threshold up to a value of $10^{-2}$. This again can be explained by the bounded actual caliber for small enough threshold values. The constructed $3$-grid method is overall more sensitive to the threshold, which is largely due to the fact that only smoothed test vectors are used. We show in~\cref{sec:bootstrap} that this problem vanishes when a more elaborate multigrid setup is used. Based on these findings, we propose to use $10^{-2}$ as a default choice for the correlation threshold.

\begin{figure}[htb]
\centerline{\includegraphics[width=1\textwidth]{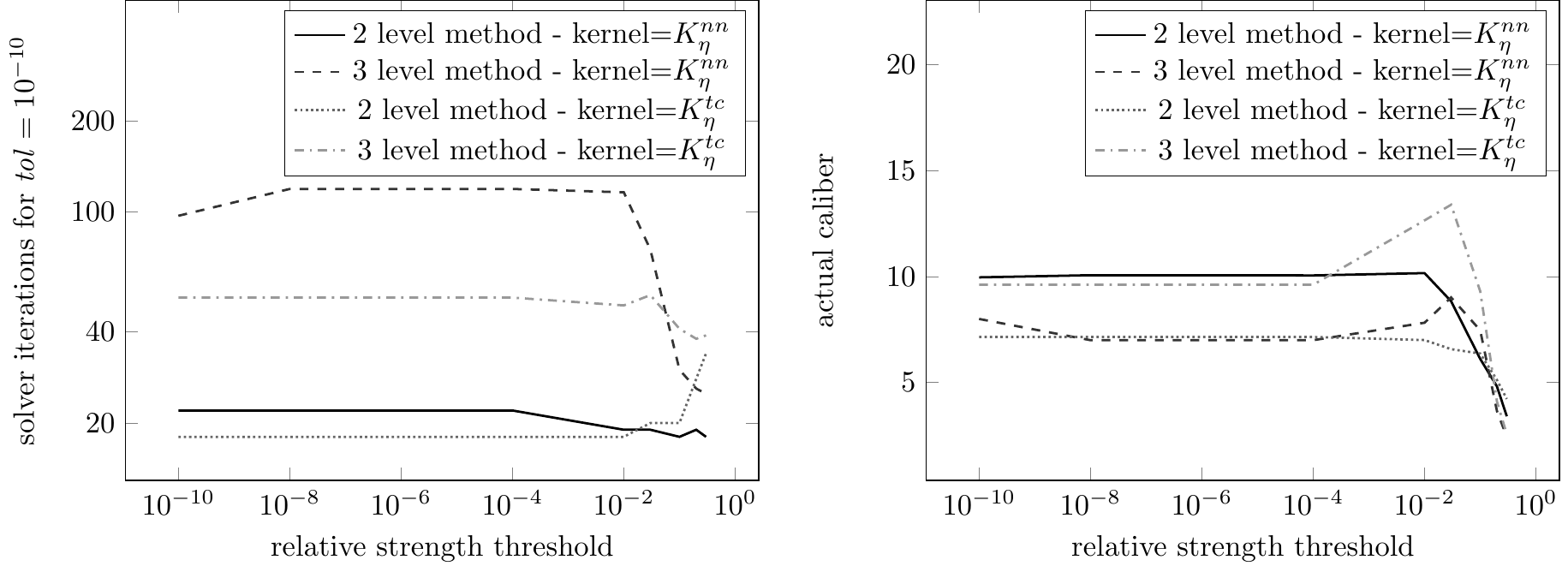}}
\caption{Study of iteration count (left) and actual caliber (right) with respect to the relative strength threshold. \label{fig:lowcut_study}}
\end{figure}
Similar results are obtained for the relative strength threshold as depicted in~\cref{fig:lowcut_study}. A large stable plateau can be used to fix the default value of this parameter to $10^{-2}$ as well.

\begin{figure}[htb]
\centerline{\includegraphics[width=1\textwidth]{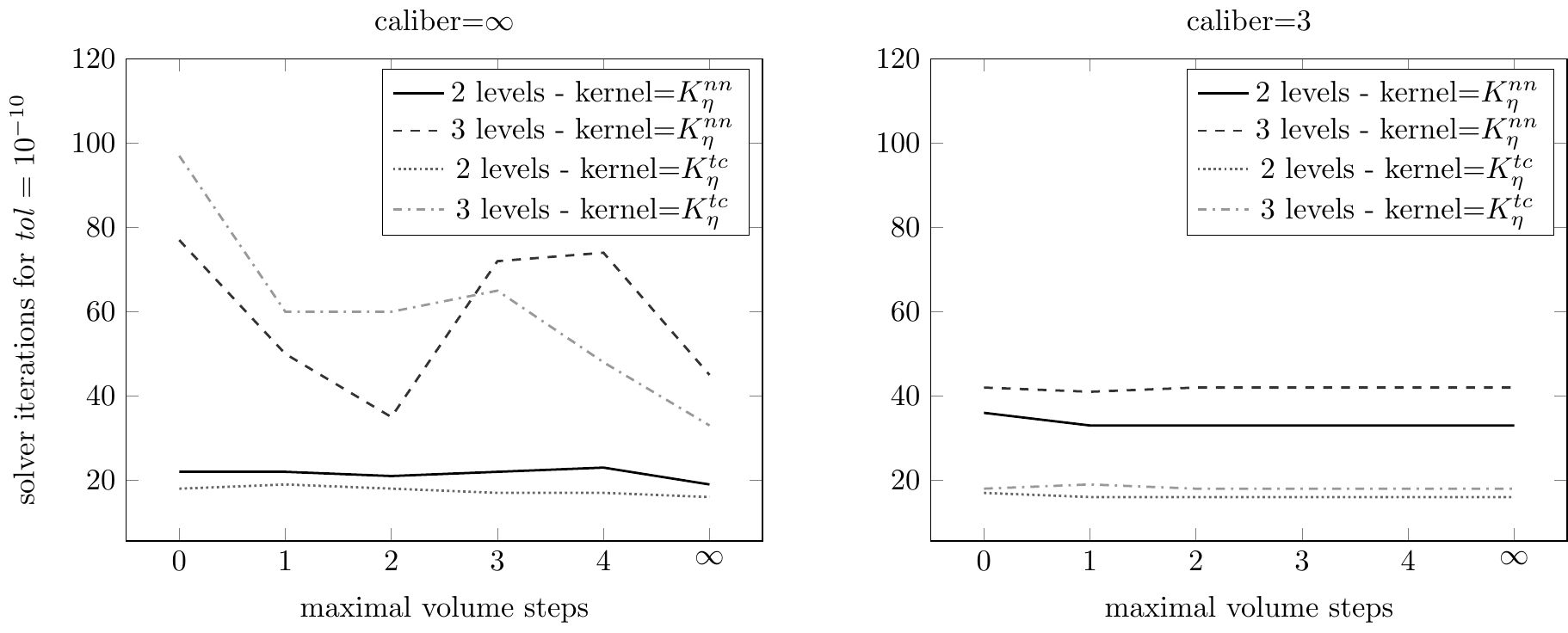}}
\caption{\emph{left:} maximal volume study, \emph{right:} relative weight low cut study.}
\label{fig:maximal_volume_study}
\end{figure}
Last, we try to gauge the influence of the maximal volume post-processing on the overall coarsening process and the quality of the corresponding 2- and 3-grid methods. In~\cref{fig:maximal_volume_study} we see that it does not affect convergence if the method converged fast before the maximal volume corrections. In case the method did not converge fast we observe overall small improvements, but it is hard to judge from this data if the post-processing by maximal volume has a significant effect. However, when inspecting the proposed coarse variable sets we do find marked visual "improvements" on the homogeneity of the obtained variable distributions. We thus propose to use the maximal volume correction and maximal number of iterations (maxvol followed by \LARS{}) to $4$. Note, that the complexity of the maximal volume correction solely depends on the number of variable swaps in and out of $\mathcal{C}$ and we typically find that only few swaps are needed for the tests considered here and in~\cref{sec:bootstrap}.

To summarize the parameter study, we collect all parameters along their default values in~\cref{tab:params}.

\begin{table}[ht]
\centerline{
\begin{tabular}{|l|l|l|l|}
\hline
                  & parameter name                &   default value  \\ \hline
kernel            & distance measure              & graph distance \\
                  & kernel function               &   $K_\eta^{tc}$  \\
                  & kernel radius                 &  $4$ [edges]    \\ \hline
\LARS{}           & correlation threshold         &  $10^{-2}$            \\
                  & rel.\ strength threshold      &  $10^{-2}$      \\
                  & caliber threshold             &  $3$      \\
                  & sign constraint               & \emph{false}   \\ \hline
coarsening        & maxvol iterations             & $4$ \\\hline \hline 
          AMG    & smoother                       &  Gauss-Seidel   \\
                 & number of test vectors         &  8              \\
                 & initial smoothing              &  4 [iterations] \\
                 & V-cycle pre-smoothing          &  1 [iterations] \\
                 & V-cycle post-smoothing         &  1 [iterations] \\
\hline
\end{tabular}
}
\caption{Summary of all parameters}
\label{tab:params}
\end{table}

\begin{figure}[htb]
\centerline{Poisson's equation, Gauss-Seidel smoothing} 
\centerline{\includegraphics[width=1\textwidth]{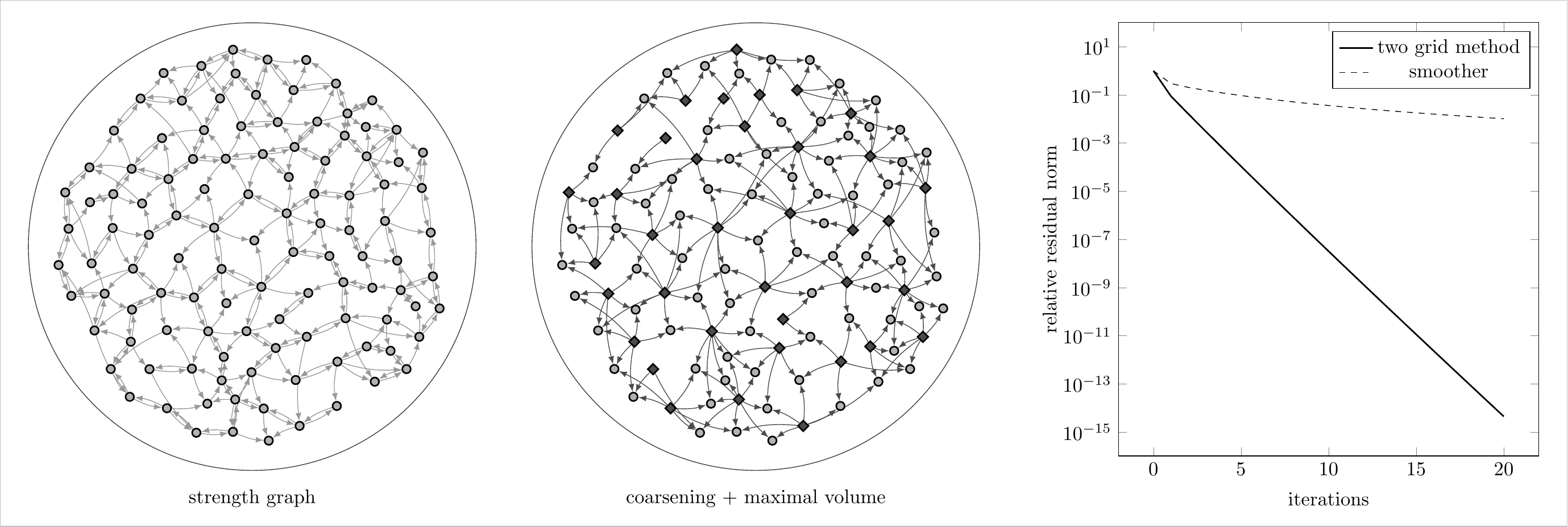}}
\vspace{.5em}
\centerline{Poisson's equation with anisotropy, Gauss-Seidel smoothing} 
\centerline{\includegraphics[width=1\textwidth]{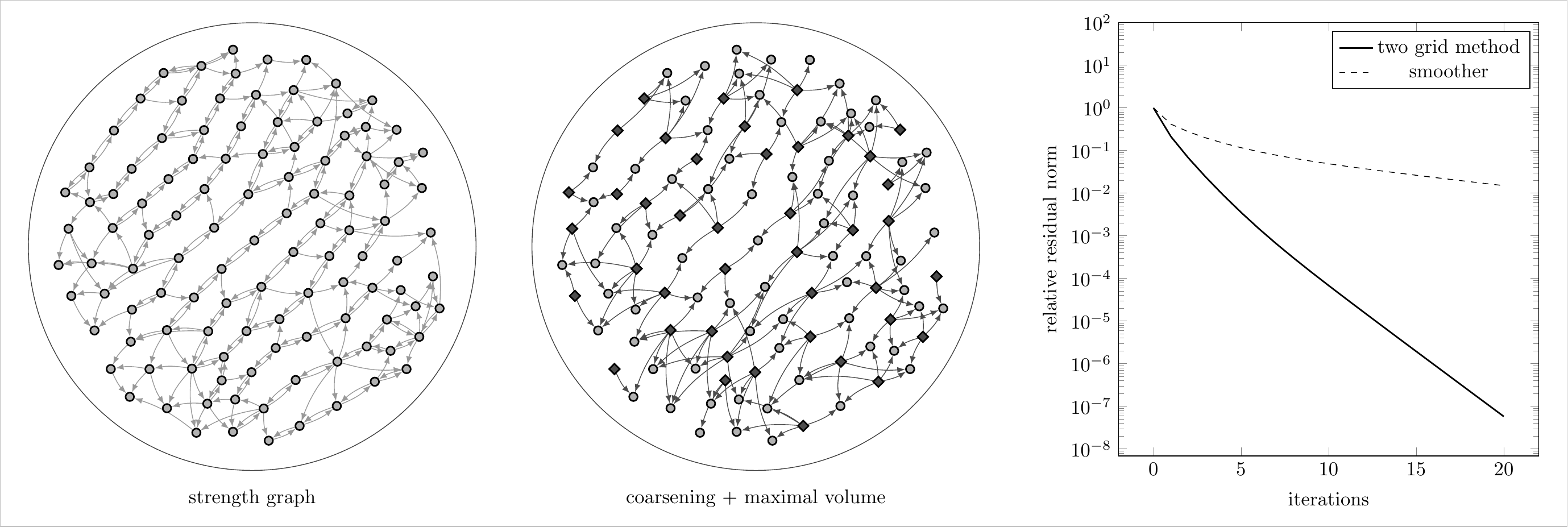}}
\vspace{.5em}
\centerline{Poisson's equation, block Gauss-Seidel smoothing} 
\centerline{\includegraphics[width=1\textwidth]{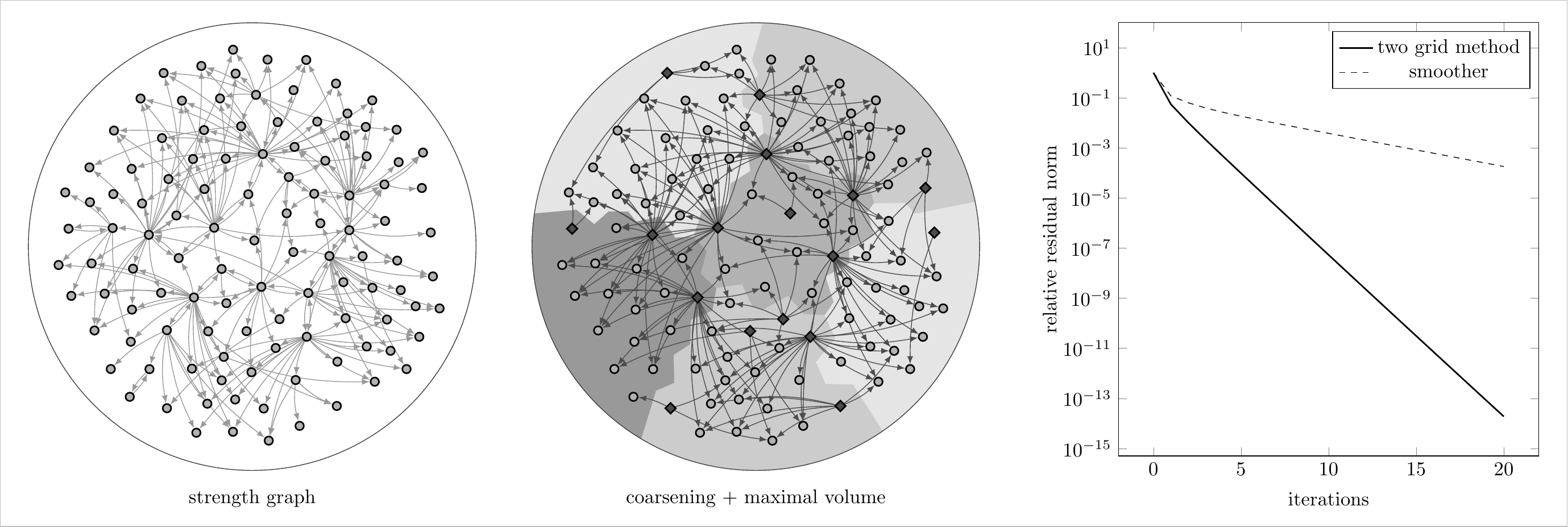}}
\caption{Illustration of the interpolation coupling structure computed by \LARS{} coarsening for three different cases. (top) Poisson w/o anisotropy, Gauss-Seidel; (middle) Poisson with anisotropy ($\alpha=\tfrac{\pi}{4}$, $\varepsilon=0.01$); (bottom) Poisson w/o anisotropy, block Gauss-Seidel (6 blocks, 4 colors)
\label{fig:poisson1}}
\end{figure}

\paragraph{Qualitative analysis}
Using the default parameters we can now take an in-depth look at the constructed strength graphs, coarse variable sets, and interpolation relations. To this end we consider an even smaller problem with $h_{\rm max}=0.2$, which yields a matrix $A\in \mathbb{R}^{112\times 112}$. This time we consider both problems with and without anisotropy as shown in~\cref{fig:poisson1}. In all tests we reduced the number of test vectors to $K=8$.

The first row of~\cref{fig:poisson1} contains results for the problem without anisotropy using Gauss-Seidel smoothing. \CLARS{} coarsening yields 37 coarse grid variables resulting in a coarsening ratio of approx.\ $0.33$ which is a reasonable choice for a caliber of three and a triangulated domain. The asymptotic convergence rate of the corresponding two-grid method is approximately $0.20$, which is not ideal, but for such a simple setup not bad either.

The second row replaces the problem by its anisotropic version with $\alpha=\frac{\pi}{4}$ and $\varepsilon=0.01$, while using Gauss-Seidel smoothing again. Clearly, \LARS{} coarsening is able to discover the anisotropy and in turn constructs a coarsening and interpolatory relations that follow it nicely. The convergence is worse than for the problem without anisotropy, as to be expected, but still far better than the standalone smoother.

Last, in the third row we tested \LARS{} coarsening for the problem without anisotropy using a block Gauss-Seidel smoother. Its blocks are depicted as the shaded regions of the grid and only one iteration of this smoother is used to generate the test vectors. For this test we chose to use the nearest-neighbor kernel as the tri-cube kernel reduces the influence of information on a scale shorter than the block-size.\footnote{A generalization of the tri-cube kernel, which uses the block geometry could be beneficial here.} It can be clearly seen that there are a small number of variables that many other variables want to interpolate from. Most of these variables are located close to the block boundaries of the smoother. The convergence rate is roughly $0.20$, as for the first row. The resulting coarsening ratio of $0.15$ yields a much coarser grid which compensates the cost of the more expensive smoother. This shows that the interpolation adapted itself well to the stronger smoother. 

In all three tests the placement of coarse grid variables and corresponding interpolatory relations becomes somewhat chaotic in proximity to the circle boundary. This does not seem to influence convergence dramatically and its effect would be further reduced when increasing the problem size.

\subsection{\CLARS{} coarsening in bootstrap AMG}\label{sec:bootstrap}
One benefit of \LARS{} coarsening is that it can be seamlessly integrated into the bootstrap algebraic multigrid framework as it uses the same test vectors that are needed for \leastsquares{} interpolation. In order to present some multigrid tests using \LARS{} coarsening in the framework, we briefly review the overall bootstrap approach, introduced in~\cite{BranBranKahlLivs2011} and modified in~\cite{BranKahlFalgHuCao2017}. As we have already introduced the concept of \leastsquares{} interpolation in~\cref{sec:lsinterpolation}, we only have to explain the bootstrap setup cycle.


The first leg of a bootstrap multigrid setup $V$-cycle consists of the method introduced in~\cref{sec:larcoarsening}. That is, based on smoothed test vectors a first multigrid hierarchy is constructed by \LARS{} coarsening and \leastsquares{} interpolation.

Now let $P_{i}^{i+1}$ denote the interpolation that maps vectors from grid $i+1$ to grid $i$ and $ P := \prod_{i=1}^{L-1} P_i^{i+1} $ the concatenation of all these interpolations which yields an interpolation operator that maps from the coarsest to the finest level. 
Based on the initial multigrid hierarchy, the bootstrap setup now considers a generalized eigenvalue problem on the coarsest grid\footnote{We chose to present the original bootstrap idea here for the sake of simplicity. According to~\cite{BranKahlFalgHuCao2017} the eigenvectors computed on the coarsest grid should be tied to the symmetrized smoothing operator $I - \widetilde{M}^{-1}A$.}
\begin{equation} \label{eq:coarsegev}
A_c V = \lambda P^T P V.
\end{equation} The eigenvectors to small eigenvalues of this eigenproblem are then interpolated through the multigrid hierarchy, where smoothing is applied to them on every grid, to augment the set of test vectors. In case additional setup iterations should be carried out, these vectors are combined with the smoothed random test vectors, which are then used to build the next (improved) multigrid hierarchy. For details of the bootstrap setup we refer to~\cite{BranBranKahlLivs2011}.



For the problem size scaling test presented in~\cref{tab:AMGscaling}, the parameters we use are the default ones from~\cref{tab:params} except for the number of test vectors which we increased to $K=16$. On the coarsest grid we computed the $16$ eigenvectors corresponding to the smallest eigenvalues of~\cref{eq:coarsegev}. As a measure of efficiency of the resulting V$(1,1)$-cycle, we report the number of preconditioned CG iterations needed to reduce the initial residual norm by a factor of $10^{10}$. As we can see the method scales nicely up to $135,\!777$ unknowns and a $7$ level method as long as we adjust the number of bootstrap cycles when scaling the problem. Note that the size of the generalized eigenvalue problem on the coarsest grid remains roughly constant implying that both the setup and the solve routine preserve the optimal complexity $\mathcal{O}(n)$ per iteration for the solver.

\begin{table}[ht]
\centerline{
\begin{tabular}{|l|l|l|l|l|l|}
\hline
$h_\mathit{max}$ & system size $n$ & levels & AMG iter & setup iter &  CG iter \\ \hline
$0.2$            & $112$           & $2$    & $9$      & $1$        &  $41$    \\
$0.1$            & $480$           & $3$    & $10$     & $1$        &  $81$    \\
$0.05$           & $2,\!045$       & $4$    & $10$     & $2$        &  $152$   \\
$0.025$          & $8,\!593$       & $5$    & $10$     & $2$        &  $286$   \\
$0.0125$         & $34,\!348$      & $6$    & $11$     & $2$        &  $501$   \\
$0.00625$        & $135,\!777$     & $7$    & $11$     & $3$        &  $892$   \\
\hline
\end{tabular}
}
\caption{Scaling test of an AMG V-cycle with one step of pre- and post-smoothing. $h_\mathit{max}$ denotes the maximum edge length. Halving $h_{\rm max}$ increases the number of unknowns by a factor of roughly $4$. We report the number of multigrid preconditioned CG iterations to reduce the initial residual norm by a factor of $10^{10}$.\label{tab:AMGscaling}}
\end{table}

\begin{figure}[htb]
\centerline{\includegraphics[width=1\textwidth]{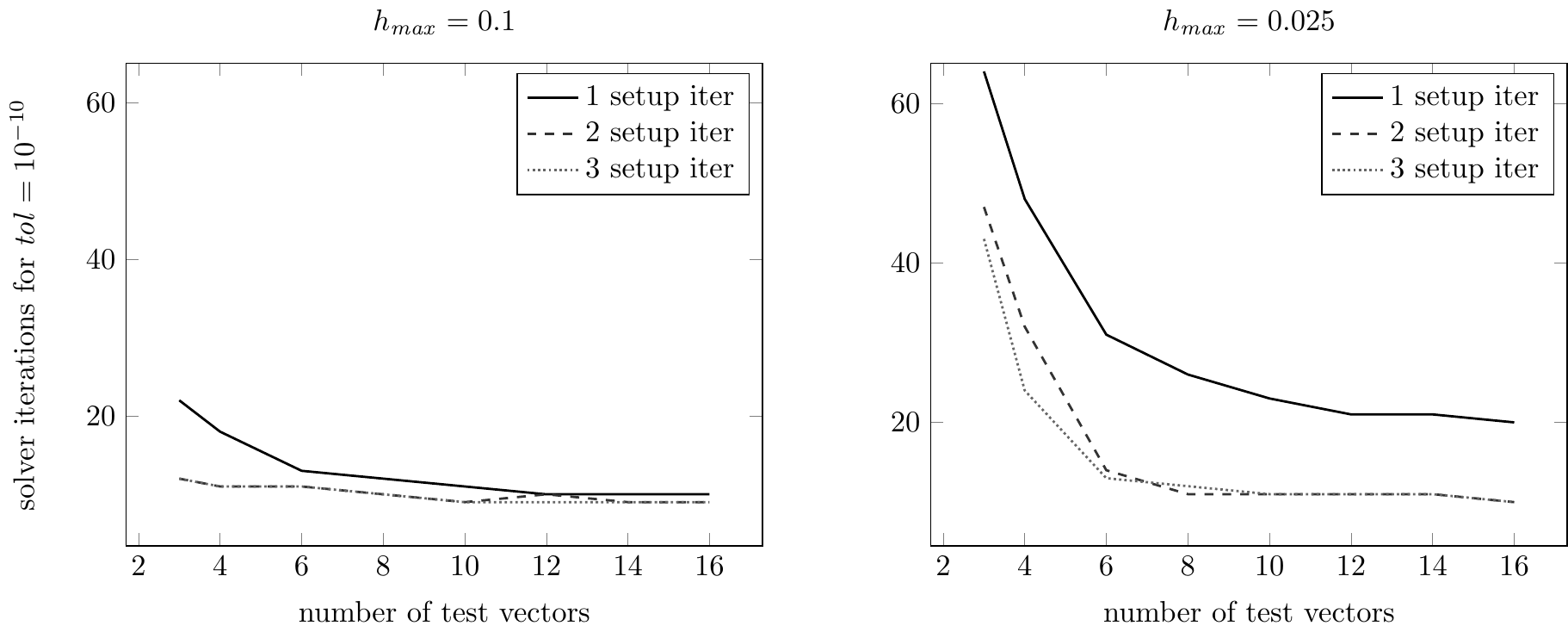}}
\caption{Preconditioned CG iterations as a function of the number of test vectors: (left) $h_\mathit{max}=0.1$; (right) $h_\mathit{max}=0.025$}
\label{fig:testvectors}
\end{figure}

These results are meant as a proof of concept. Clearly, the number of test vectors needs to be tuned in order to optimize the resulting multigrid method (setup+solve). In order to give some insight into the scaling of the method with respect to the number of test vectors used, we report multigrid preconditioned CG iteration counts for different choices of $K$ for two grid sizes in~\cref{fig:testvectors}. Clearly, the first setup cycle, which does not use any multigrid enhanced test vectors, profits most from additional test vectors. In addition there seems to be a lower bound on the number of test vectors needed to establish a stable plateau at around $K=6$. This is not surprising as at a caliber of $3$, fewer than $6$ test vectors might lead to severe overfitting.


\section{Conclusion}

We have shown as a proof of concept that the coarsening problem of algebraic multigrid can be solved using a machine learning ansatz. Our \LARS{} coarsening approach demonstrates that it is possible to find meaningful and efficient coarsenings without the need to deviate from the bootstrap algebraic multigrid framework and that it can be integrated seamlessly. That is, it can work with few test vectors and few setup smoothing iterations to determine a local model of algebraically smooth error. In this it is able to determine interpolation sets $\mathcal{C}_{i}$ in an almost greedy fashion guided by the $\ell_{1}$ constraint \leastsquares{} problem. It thus circumvents some of the problems inherent in approaches that consider only binary variable relations, where interactions between interpolating variables are neglected, without the need to consider all possible candidate subsets.

The proposed method not only showed that it yields scalable methods for constant coefficient problems on unstructured meshes, but also that it is able to generate suitable coarsenings for anisotropic problems and for more complex smoothers such as block smoothers, where interpolation relations are much more difficult, if at all possible, to describe geometrically.
As the method is based solely on local operations it is suitable for parallelization and we expect to make further progress and improvements by considering it for systems of partial differential equations and overlapping block smoothers.

Our results indicate that a machine learning perspective of algebraic multigrid, that views the construction of interpolation as a problem of learning the local nature of algebraically smooth errors, is helpful for future developments of adaptive algebraic multigrid methods.

\subsection*{Acknowledgements}

We would like to thank J.~Brannick and A.~Frommer for discussions and helpful remarks in preparing this manuscript\textbf{}.

\printbibliography


\typeout{get arXiv to do 4 passes: Label(s) may have changed. Rerun}

\end{document}